\documentclass{article}
\usepackage{graphicx}
 \usepackage{mathptmx}
\usepackage{amsmath, amstext,amssymb,amsfonts}
\usepackage[english]{babel}
\usepackage{delarray}
\usepackage{mathptmx}
\usepackage{amsmath, amstext,amssymb,amsfonts}
\usepackage[english]{babel}
\usepackage{delarray}

\DeclareFontFamily{U}{mathx}{\hyphenchar\font45}
\DeclareFontShape{U}{mathx}{m}{n}{
      <5> <6> <7> <8> <9> <10>
      <10.95> <12> <14.4> <17.28> <20.74> <24.88>
      mathx10
      }{}
\DeclareSymbolFont{mathx}{U}{mathx}{m}{n}
\DeclareFontSubstitution{U}{mathx}{m}{n}
\DeclareMathAccent{\widecheck}{0}{mathx}{"71}
\DeclareMathAccent{\wideparen}{0}{mathx}{"75}

\numberwithin{equation}{section}
  \newtheorem{teo}{Theorem}
\newtheorem{corr}[teo]{Corollary}

 \newcommand{\cb}{CB(\mathbb{R})}
  \newcommand{\rpl}{\mathbb{R}^2_{+}}  \newcommand{\cpl}{\mathbb{C}_{+}}
   \newcommand{\R}{\mathbb R}    \newcommand{\Co}{\mathbb C}

\date{}
\title{A note on harmonic continuation of characteristic function}
\author{\large{ Saulius  Norvidas}}
\date{\footnotesize Institute of Data Science and Digital Technologies, Vilnius University, \\ Akademijos str. 4, Vilnius LT-04812, Lithuania\\
 ({\rm{e-mail: norvidas{@}gmail.com}})}
\begin{document}

\maketitle
 {{ {\bf Abstract}}}
We propose a necessary and sufficient condition for a real-valued function on the real line to be a characteristic function of a probability measures.  The statement is given  in terms of harmonic  functions and completely monotonic functions.



\section{ Introduction }

{\large{

If $\sigma$ is a probability measure on the real line $\mathbb{R}$, then
\[
\theta(x)=\int_{-\infty}^{\infty}e^{ixt}\,d\sigma(t)
\]
is called the characteristic function of the probability measure $\sigma$. We will call $\theta$ the characteristic function for short. Note that a characteristic function  is real-valued if and only if  it is even  on $\R$.

\ \ \ \ The motivation of the present paper stems from Egorov's paper  \cite[p. 567]{2}, where  the following criterion for real-valued infinitely differentiable characteristic  functions is given:

\begin{teo}\label{E}. Let $\varphi(t)$ be  an infinitely differentiable absolutely integrable even
function  on $\R$ with real values. Let  $\varphi(t)$ satisfy the condition $\varphi^{(m)}(t) t^{-2}={O(1)}$ as $t\to\infty$, $m=1,2,\dots$,  and let the Fourier transform $\widehat{\varphi}$  belong to $L^1(\mathbb{R})$. Then $\varphi(t)$  is a characteristic function if and only if

\ \ \ 1) \ for all $0 < \delta < \infty$ and all $p = 0, 1, 2, . . .$
\[
(-1)^p\int_0^{\infty}\frac1{\delta+t^2}\varphi^{(2p)}(t)\,dt\ge 0,
\]
\[
(-1)^{p+1}\int_0^{\infty}\frac{t}{\delta+t^2}\varphi^{(2p+1)}(t)\,dt\ge 0;
\]

\ \ \ 2)\ $\varphi(0)=1$.
\end{teo}

\ \ \ \ We claim that some conditions of Theorem 1 are unnecessary. Furthermore, this criterion can not be applied to the simplest infinitely differentiable characteristic functions  $\varphi(t)\equiv 1$ and $\varphi(t)=\cos at$, where $a\in\R$ (since they  are not  absolutely integrable on $\R$).  Therefore, we will follow here a  different approach to a similar criterion associated with the theory of harmonic functions. This approach allow us to simplify and strengthen the statement of Theorem 1; compare Corollary 3  below.

\ \ \ \ Let us start by introducing some  notation and  basic facts  that will be used
throughout this paper. Let $CB(\mathbb{R})$ denote the Banach space of bounded continuous  functions $f: \mathbb{R}\to\mathbb{R}$ with the usual uniform norm. If $f\in\cb$, then, as is well known (see, for example \cite[Chapter II]{5}), the Dirichlet problem for the upper half plane $\rpl=\{ (x,y)\in \mathbb{R}^2:\ y>0\}$ has  only unique solution  in $CB(\overline{\rpl})$, where  $\overline{\rpl}$  is the closure of $\rpl$. More precisely, given $f\in\cb$, there exists a  function $u_f$  on  $\overline{\rpl}$ such that $u_f$ is bounded and  continuous in $\overline{\rpl}$,
\begin{equation}\label{1.1}
\lim_{y\to 0} u_f(x,y)=f(x)
\end{equation}
for each  $x\in\mathbb{R}$, and $u_f$  is harmonic  in $\rpl$. This means that  $u_f$ is infinitely differentiable in $\rpl$ and there  satisfies
\[
\biggl(\frac{\partial^2}{\partial x^2}+\frac{\partial^2}{\partial y^2}\biggr)\, u_f=0.
\]
 We call $u_f$ the harmonic continuation of $f$ into $\rpl$. The solution of the Dirichlet problem for $\rpl$ can be obtained  by using  the Poisson kernel
\begin{equation}\label{1.2}
P(x,y)=\frac{1}{\pi}\frac{y}{x^2+y^2}.
\end{equation}
 Namely, the harmonic continuation $u_f$ can be obtained  as a  convolution with respect to $x$ of $f$ and $P(x,y)$
\begin{equation}\label{1.3}
u_f(x,y)=\Bigl(f\ast_{x} P\Bigr)(x,y)=\frac1{\pi}\int^{\infty}_{-\infty}\frac{y}{(x-t)^2+y^2}f(t)\,dt.
\end{equation}
where $(x,y)\in \rpl$.

\ \ \ \ We recall that a function $\omega: [0, \infty)\to \R$ is called completely monotonic if it is infinitely differentiable on $[0,\infty)$ and
\begin{equation}\label{1.4}
(-1)^n \omega^{(n)}(x)\ge 0
\end{equation}
for $x\in(0, \infty)$ and all  $n=0,1,2,\dots$. The following elementary functions are immediate examples of completely monotonic functions, which
is verified directly:
\[
e^{-ax},\qquad \frac1{(\alpha x+\beta)^p},\quad {\text{and}}\quad \ln\Bigl(b +\frac{c}{x}\Bigr),
\]
where $a\ge 0$, $\alpha\ge 0$, $\beta\ge 0$,  and $p\ge 0$  with $\alpha$ and $\beta$ not both zero and $b\ge 1$, $c > 0$.

\begin{teo}\label{t2}.
Suppose that $f\in \cb$ is an even function  and  $f(0)=1$. Then $f$ is a characteristic function if and only if the restriction  of  $u_f$ to the imaginary axis, i.e., the function  $y\to u_f(0,y)$, $y\in [0,\infty)$, is completely monotonic.
\end{teo}

\begin{corr}\label{c3}.
Let $f\in \cb$. Suppose

(i)\  $f$ is  infinitely differentiable  on $\R$ and $ f^{(k)}\in\cb$ for all $k=1,2\dots$;

(ii)\  $f$ is even  and  $f(0)=1$.

Then $f$ is characteristic function if and only if
\begin{equation}\label{1.5}
\int_{-\infty}^{\infty}\Im\Bigl(\frac1{t-iy}\Bigr)\,f(t)\,dt\ge 0,
\end{equation}
\begin{equation}\label{1.6}
(-1)^{n+1}\int_{-\infty}^{\infty}\Re\Bigl(\frac1{(t-iy)^2}\Bigr)\,f^{(2n)}(t)\,dt\ge 0
\end{equation}
and
\begin{equation}\label{1.7}
(-1)^{n}\int_{-\infty}^{\infty}\Im\Bigl(\frac1{(t-iy)^2}\Bigr)\,f^{(2n+1)}(t)\,dt\ge 0
\end{equation}
for   $y>0$ and all $n=0,1,2,\dots$.
\end{corr}

}

\section{ PRELIMINARIES AND PROOFS}\label{s:2}
\large{
\ \ \ \ A complex-valued function $\varphi$ on $\R$ is said to be positive definite if
\[
\sum_{i,j=1}^{n} \varphi(x_i-x_j) c_i\overline{c_j}\ge 0
\]
for every choice of $x_1,\dots,,x_n\in\R$, for every choice of complex numbers $c_1,\dots, c_n\in\Co$, and all $n\in\mathbb{N}$. The Bochner theorem (see \cite[p. 150]{3}) characterizes continuous positive definite functions: a continuous function $\varphi: \R\to \Co$ is positive definite if and only if there exists  a finite non-negative measure $\mu$ on $\R$ such that
 \[
\varphi(x)=\int_{-\infty}^{\infty} e^{ixt}\,d\mu(t).
\]
Note that  $\varphi$ is a characteristic function if and only if $\varphi$ is  continuous positive definite  and $\varphi(0)=1$. Any characteristic function $\varphi$ satisfies:

\ \ (i)\ $\varphi$ is uniformly continuous on $\R$;

\ \ (ii)\  $\varphi$ is bounded  on $\R$, i.e.,  $\varphi(x)|\le \varphi(0)=1$ for all $x\in\R$;

\ \ (iii)\ The real part of  $\theta(x)$ is also an even characteristic function.

\ \ \ \ In the case of completely monotonic functions, the Bernstein-Widder theorem (see \cite[p.  161]{7}) asserts that a function $f: [0,\infty)\to\R$ is completely monotonic if and only if it is the Laplace transform of a finite non-negative measure $\eta$ supported on $[0,\infty)$, i.e.,
\begin{equation}\label{2.1}
\varphi(x)=\int_{0}^{\infty} e^{-xt}\,d\eta(t)
\end{equation}
for $x\in [0,\infty)$.

\ \ \ \ Let $f\in\cb$. It is well known that  harmonic continuation (1.3) has harmonic conjugate. Recall that if $U$ is a harmonic function in $\rpl$, then another  function $V$ harmonic in $\rpl$ is called  harmonic conjugate of $U$, provided  $U+iV$ is analytic in  $\cpl=\{z=x+iy\in\mathbb{C}:\ y>0\}$. The harmonic conjugate $V$ is unique, up to adding a constant. Let $v_f$ denote a harmonic conjugate of $u_f$.  We can choose (see  \cite[p.p. 108-109]{4})
\begin{equation}\label{2.2}
v_f(x,y)= \frac1{\pi}\int_{-\infty}^{\infty} \Bigl(\frac{x-t}{(x-t)^2+y^2}+\frac{t}{t^2+1}\Bigr) f(t)\,dt,
\end{equation}
where $x\in\R$ and $y>0$. Then
\begin{equation}\label{2.3}
E_f(z)=E_f(x,y)=u _f(x,y)+ iv_f(x,y)
\end{equation}
 is analytic in $\cpl$.

\ \ \ \  We note that usually the following simpler integral
\begin{equation}\label{2.4}
\widetilde{v}_f(x,y)= \frac1{\pi}\int_{-\infty}^{\infty} \frac{x-t}{(x-t)^2+y^2} f(t)\,dt
\end{equation}
chosen as the harmonic conjugate of $u_f$.  In that case  the analytic function $u_f+i \widetilde{v}_f$ coincides  with the usual Cauchy type integral
\[
(u_f+i\widetilde{v}_f)(z) = \frac{i}{\pi}\int_{-\infty}^{\infty} \frac{f(t)}{z-t}\,dt
\]
for  $z\in \cpl$. However, the integral in (\ref{2.4}) converges as long as
\[
\int_{-\infty}^{\infty} \frac{|f(t)|}{1+|t|}\,dt<\infty.
\]
 For example, this condition is satisfied if $f\in L^p(\R)$, $1 \le p<\infty$, but not in the case of an arbitrary  $f\in \cb$.

\ \ \ {\bf Proof  of Theorem \ref{t2}}. \  Suppose that  $f$ is a real-valued  characteristic function. By the Bochner  theorem, there exists a probability measure $\sigma$ on $\R$ such that
\begin{equation}\label{2.5}
f(t)=\int_{-\infty}^{\infty} e^{ixt}\,d\sigma(x).
\end{equation}
For each $y>0$ the Poisson  kernel (\ref{1.2}) is a Lebesgue integrable function of $x$ on   $\R$. Hence, by the Fubini theorem, it follows from (\ref{1.3}) and (\ref{2.5}) that
\[
u_f(0,y)=\frac1{\pi}\int_{-\infty}^{\infty}\Bigl(\int_{-\infty}^{\infty} P(t,y)e^{ixt}\,dt\Bigr)\,d\sigma(x).
\]
Since
\[
\int_{-\infty}^{\infty} P(t,y)e^{ixt}\,dt=e^{-|x|y},
\]
it follows that
\begin{equation}\label{2.6}
u_f(0,y)=\int_{-\infty}^{\infty} e^{-|x|y}\,d\sigma(x).
\end{equation}
For each $y>0$ and any $n=0,1,1,2\dots$, there exists $0<a(y, n)<\infty$ such that
\[
\sup_{x\in\R}\ \frac{\partial^{n}}{\partial y^n}e^{-|x|y}\le a(y,n).
\]
Therefore, if  we recall that  $\sigma$ is a finite measure, then we have that  the partial derivatives of $u_f(0,y)$ in $y$  can be obtained by differentiation under the
integral sign in (\ref{2.6}) (see, for example,  \cite[p. 283]{1}).  Now the positiveness $\sigma$ implies that
\[
(-1)^{n}\Bigl(\frac{\partial^n}{\partial y^n}\, u_f\Bigr)(0,y)=\int_{-\infty}^{\infty}|x|^n e^{-|x|y}\,d\sigma(x)\ge 0
\]
for all $y>0$ and all $n=0,1,2,\dots$. Note that the harmonic function $u_f(x,y)$ is continuous in $\overline{\rpl}$. Therefore, by the Bernstein-Widder theorem,  we obtain  that  $u_f(0,y)$ is completely monotonic on $[0,\infty)$.

\ \ \ \ Suppose $u_f(0,y)$ is completely monotonic for $y\in [0,\infty)$. By the Bernstein-Widder theorem, there exists a finite non-negative measure $\eta$ supported on $[0,\infty)$ such that
\begin{equation}\label{2.7}
u_f(0,y)=\int_{0}^{\infty} e^{-yt}\,d\eta(t),
\end{equation}
$y\in[0,\infty)$. According to  (\ref{1.1}), we get
\begin{equation}\label{2.8}
\lim_{y\to 0} u_f(0,y)= f(0)=1.
\end{equation}
This means that $\eta$ is a probability measure. Set
\begin{equation}\label{2.9}
K_{\eta}(z)=\int_0^{\infty} e^{izt}\,d\eta(t)
\end{equation}
for  $z\in \overline{\mathbb{C}}_{+}=\{ z=x+iy\in \mathbb{C}: \ y\ge 0\}$. We can consider the measure $\eta$ as tempered distribution on $\R$. Then $K_{\eta}$ is the distributional Laplace transform of $\eta$ (see \cite[p. 127]{6}).  Therefore, since $\eta$ is supported on $[0,\infty)$, we have that $K_{\eta}$ is analytic function in $\R+i(0,\infty)={\mathbb{C}}_{+}$ and the following differentation formula in  ${\mathbb{C}}_{+}$ holds
\[
\frac{d^n}{d\,z^n}  K_{\eta}(z)=(i)^n\int_0^{\infty} e^{ixt}\Bigl( t^ne^{-yt}\Bigr)\,d\eta(t),
\]
 \cite[p.p. 127-128]{6}.

\ \ \ \ We claim that $K_{\eta}$ and  (\ref{2.3}) coincide in $\cpl$. Indeed, since $f$ is even on $\R$, it follows from (\ref{2.2}) that $v_f(0,y)=0$ for each $y>0$. Hence
\begin{equation}\label{2.10}
E_f(iy)= u_f(0,y)
\end{equation}
for $y>0$. On the other hand, using  (\ref{2.7}) and (\ref{2.9}), we have
\[
K_{\eta}(iy)= u_f(0,y),
\]
 $y> 0$. If we combine this with   (\ref{2.10}) and use the uniqueness theorem for analytic functions, we obtain that $E_f=K_{\eta}$ in $\cpl$.

\ \ \ \ By   (\ref{2.9}), we have that for any fixed $y_0\in(0,\infty)$ the function
$x\to K_{\eta}(x+iy_0)$
is   positive definite  for $x\in\R$. Therefore the function $x\to \Re \bigl[K_{\eta}(x+iy_0)\bigr]$ also is  positive definite  for $x\in\R$. Since $E_f=K_{\eta}$ in $\cpl$, it follows from (\ref{2.3}) that
\[
u_f(x,y_0)= \Re \bigl[E_f(x,y_0)\bigr]=\Re \bigl[K_{\eta}(x+iy_0)\bigr]
\]
for $x\in\R$. Thus each function $x\to u_f(x,y_0)$, where $y_0\in(0,\infty)$, is  positive definite   on $\R$.  Since  the pointwise limit of positive definite functions also is positive definite  function,    (\ref{1.1}) implies that  $f$ is  positive definite  on $\R$. Finally, since $f$ is continuous on $\R$, it follows from (\ref{2.8}) and the Bochner theorem that $f$ is a characteristic function. This proves Theorem 2.

\ \ \ {\bf Proof  of  Corollary  \ref{c3}}. \ By Theorem 2 and (1.4), a necessary and sufficient condition for $f$ to be a
 characteristic function is that
\begin{equation}\label{2.11}
(-1)^k\Bigl[\frac{\partial^k}{\partial y^k} u_f(0,y)\Bigr]\ge 0
\end{equation}
for $k=0,1,2,\dots$, and $y>0$. In the case $k=0$, the condition (2.11) is
\begin{equation}\label{2.12}
u_f(0,y)=\frac1{\pi}\int_{-\infty}^{\infty}\frac{y}{t^2+y^2} f(t)\,dt=\frac1{\pi}\int_{-\infty}^{\infty}\Im\Bigl(\frac1{t-iy}\Bigr) f(t)\,dt\ge 0.
\end{equation}

\ \ \ \ Let us calculate the partial derivatives in (\ref{2.11}). To this end, we recall that if $F$ is a analytic function, then  the Cauchy-Riemann equations  imply
\[
\frac{d}{dz} F(z)=\frac{\partial}{\partial y}\Bigl( \Im F(z)\Bigr)-i\frac{\partial}{\partial y}\Bigl( \Re F(z)\Bigr).
\]
Hence, in the case if $F$ is  the  function (\ref{2.3}), then we get
\[
\frac{d^{2n}}{dz^{2n}} E_f(z)=(-1)^n \Bigl[ \frac{\partial^{2n}}{\partial y^{2n}}\, u_f(x,y)+ i\frac{\partial^{2n}}{\partial y^{2n}}\, v_f(x,y)\Bigr]
\]
for $n=1,2,\dots$,  and
\[
\frac{d^{2n+1}}{dz^{2n+1}} E_f(z)=(-1)^n \Bigl[ \frac{\partial^{2n+1}}{\partial y^{2n+1}} \, v_f(x,y) - i\frac{\partial^{2n+1}}{\partial y^{2n+1}} \, u_f(x,y)\Bigr]
\]
for $n=0,1,2,\dots$, where $z\in\cpl$. Therefore
\begin{equation}\label{2.13}
\frac{\partial^{2n}}{\partial y^{2n}}\, u_f(x,y)= (-1)^n \Re\Bigl[\frac{d^{2n}}{dz^{2n}}\ E_f(z)\Bigr]
\end{equation}
for $n=1,2,\dots$,  and
\begin{equation}\label{2.14}
\frac{\partial^{2n+1}}{\partial y^{2n+1}}\, u_f(x,y)= (-1)^{n+1} \Im\Bigl[\frac{d^{2n+1}}{dz^{2n+1}}\ E_f(z)\Bigr]
\end{equation}
for $n=0,1,2,\dots$.

\ \ \ \ On the other hand, by (\ref{1.3}), (\ref{2.2}) and (\ref{2.3}), we have
\[
E_f(z)=\frac{i}{\pi}\int_{-\infty}^{\infty}\Bigl(\frac1{z-t}+\frac{t}{t^2+1}\Bigr) f(t)\,dt.
\]
Hence
\[
\frac{d^{k}}{dz^{k}}\, E_f(z)=(-1)^k\frac{ i\, k!}{\pi}\int_{-\infty}^{\infty}\frac{f(t)}{(z-t)^{k+1}}\,dt=(-1)^k
\frac{ i }{\pi}\int_{-\infty}^{\infty}\frac{f^{(k-1)}(t)}{(z-t)^{2}}\,dt,
\]
where $k=1,2,\dots$. Therefore, (\ref{2.13}) and (\ref{2.14}) imply
\begin{gather}\label{2.15}
\frac{\partial^{2n}}{\partial y^{2n}}\, u_f(0,y)=(-1)^n\ \Re\Bigl[\frac{d^{2n}}{dz^{2n}}\, E_f(z)\Bigr](z=iy)=
(-1)^n\Re\Bigl[\frac{i}{\pi}\int_{-\infty}^{\infty}\frac{f^{(2n-1)}(t)}{(iy-t)^{2}}\,dt\Bigr]
 \nonumber      \\
 =\frac{(-1)^{n+1}}{\pi}\int_{-\infty}^{\infty}\Im\Bigl[\frac{1}{(t-iy)^2}\Bigr]{f^{(2n-1)}(t)}\,dt
\end{gather}
for $n=1,2,\dots$, and
\begin{gather}\label{2.16}
\frac{\partial^{2n+1}}{\partial y^{2n+1}}\, u_f(0,y)=(-1)^{n+1}\ \Im\Bigl[\frac{d^{2n+1}}{dz^{2n+1}}\, E_f(z)\Bigr](z=iy)
 \nonumber      \\
 =
(-1)^{n+1}\ \Im\Bigl[\frac{-i}{\pi}\int_{-\infty}^{\infty}\frac{f^{(2n)}(t)}{(iy-t)^{2}}\ dt\Bigr]
 =\frac{(-1)^{n}}{\pi}\int_{-\infty}^{\infty}\Re\Bigl[\frac{1}{(t-iy)^2}\Bigr]{f^{(2n)}(t)}\,dt
\end{gather}
where $n=0,1,2,\dots$.

 Finally, (\ref{2.12}), (2.15), and (2.16) show that (\ref{2.11}) is equivalent to the conditions
(\ref{1.5})--(\ref{1.7}).

}}
\end{document}